\documentclass[12pt,leqno,twoside]{amsart}
\usepackage{amssymb}
\usepackage{latexsym}
\usepackage{verbatim}
\usepackage{epsfig}
\usepackage{changebar}

\newtheorem{theorem}{Theorem}[section]

\newtheorem{lemma}[theorem]{Lemma}

\theoremstyle{definition}
\newtheorem{definition}[theorem]{Definition}
\theoremstyle{remark}

\newcommand{\reg}{{\rm{reg}}}

\newcommand{\proj}{{\rm{proj}}}

\newcommand{\Eu}{{\rm{Eu}}}

\newcommand{\cl}{{\rm{closure}}}
\newcommand{\grad}{\mathop{\rm{grad}}}
\newcommand{\ity}{{\infty}}
\renewcommand{\d}{{\rm{d}}}
\newcommand{\e}{\varepsilon}

\newcommand{\cW}{{\mathcal W}}

\newcommand{\bR}{{\mathbb R}}
\newcommand{\bC}{{\mathbb C}}
\newcommand{\bP}{{\mathbb P}}

\newcommand{\bv}{{\bf{v}}}
\newcommand{\bw}{{\bf{w}}}

\begin{document}

\title[Global Euler obstruction and polar invariants]
{Global Euler obstruction and polar invariants}

\author{\sc Jos\'e Seade}

\address{J.S. and A.V.: Instituto de Matem\'aticas, Unidad Cuernavaca, Universidad Nacional 
Aut\'onoma de M\'exico, Apartado postal 273-3, C.P. 62210, Cuernavaca, Morelos, M\'exico.}

\email{jseade@matem.unam.mx}\email{alberto@matcuer.unam.mx}

\author{\sc Mihai Tib\u ar}

\address{M.T.:  Math\' ematiques, UMR 8524 CNRS,
Universit\'e de Lille 1, \  59655 Villeneuve d'Ascq, France.}

\email{tibar@agat.univ-lille1.fr}

\author{\sc Alberto Verjovsky}

\thanks{Partially supported by CNRS-CONACYT (12409) Cooperation Program. 
The first and third named authors partially supported by
CONACYT grant G36357-E and DGPA (UNAM) grant IN 101 401}






\begin{abstract}
Let $Y\subset \bC^N$ be a purely dimensional, complex algebraic singular space. We define a global 
Euler obstruction $\Eu(Y)$ which extends the Euler-Poincar\'e characteristic in case of a nonsingular $Y$. Using Lefschetz pencils, we express $\Eu(Y)$ as alternating sum of global polar invariants.

\end{abstract}
%
\maketitle

\setcounter{section}{0}
\section{Introduction}\label{intro}

The local Euler obstruction, introduced by R. MacPherson in \cite{MP},
is an invariant of complex analytic varieties that plays an
essential role for studying the Chern classes of singular
varieties. It appeared in the construction of MacPherson's cycle,
as the constructible function $\Eu_X$. 
Roughly speaking, if $(X,x_0)$ is a possibly singular germ,
the local Euler obstruction of $X$ at $x_0$, $\Eu_X(x_0)$, is the
obstruction for extending a continuous stratified radial vector field around
$x_0$ in $X$ to a non-zero section of the Nash bundle over the Nash blow up $\tilde X$ of $X$ \cite{MP,Du1,BS,LT}.  

 We define in this paper a global Euler obstruction
$\Eu(Y)$ for an affine singular variety $Y\subset \bC^N$ of pure dimension $d$, in a similar manner, i.e. as the 
obstruction to extend a radial vector field, defined on the link at infinity
of $Y$, to a non-zero section of the Nash bundle. In case $Y$ is non-singular, $\Eu(Y)$ equals the Euler-Poincar\'e characteristic $\chi(Y)$.
Our obstruction can
be actually regarded as the top dimensional Chern-Mather class of $Y$.

We show here that $\Eu(Y)$ can be expressed as an alternating sum:

\begin{equation}\label{eq:main}
\Eu(Y) =  (-1)^{d} \alpha_Y^{(1)} + \cdots - \alpha_Y^{(d)} + \alpha_Y^{(d+1)},
\end{equation}
%
where the invariants $\alpha_Y^{(i)}$ are defined as follows:
 $\alpha_Y^{(1)}$ is the number of Morse points on the regular part $Y_\reg$
of a Lefschetz pencil on $Y$ and the following ones are similar numbers
defined on succesive generic hyperplane slices of $Y$. These invariants can be viewed as global polar multiplicities.
Formula (\ref{eq:main}) may therefore look analogous to L\^e-Teissier's one for the local Euler obstruction \cite{LT}, in which local polar multiplicities enter. 
Our proof of (\ref{eq:main}) has nevertheless different flavour. It relies on the repeated use of the Lefschetz  method of slicing by pencils and on a construction, inspired from the local setting of \cite{BLS, BMPS}, which is, roughly, as follows:  start from a stratified radial-at-infinity vector field and extend it first over a generic hyperplane section, then find the obstruction to extend it further. 

  We shall explain in \S \ref{s:remarks} how the known local formulas \cite{LT, BLS}, as well as our global formulas, 
 can be viewed as consequences, in appropriate settings, of the Lefschetz principle and of Dubson's definition of Euler obstruction relative to a non-characteristic open set \cite{Du1, Du2}.

 We thank J. Sch\"urmann for his comments on a preliminary version of
 this paper. He informed us that his forthcoming book \cite{Sch-book} contains a general viewpoint on Euler characteristics of constructible sheaves,
 generalizing in particular Dubson's index formula \cite{Du2}.

\section{Global Euler obstruction}\label{main}
 
 Let $Y\subset \bC^N$ be an algebraic space of pure
 dimension $d$. We consider the analytic closure $\bar Y$ of $Y$ in
 the complex projective space $\bP^N$ and denote by $H_\ity$ the
 hyperplane at infinity of the embedding $\bC^N\subset \bP^N$.  One
 may endow $\bar Y$ with a semi-analytic Whitney stratification $\cW$ such that the
 part ``at infinity'' $\bar Y\cap H_\ity$ is a union of strata. Since $\bar
 Y$ is projective and since the stratification $\cW$ is locally
 finite, it follows that $\cW$ has finitely many strata. The regular part $Y_\reg$ is the stratum (possibly not connected) of highest dimension $d$.
 
 Let $\widetilde Y$ denote the Nash blow-up of $Y$, that is:
 \[ \widetilde
 Y = \cl \{ (x,T_x Y_\reg)\mid x\in Y_\reg\}\subset Y\times G(d,N),\]
  where $G(d,N)$ is the Grassmannian of complex $d$-planes in $\bC^N $.
 Let $\nu :\widetilde Y \to Y$ denote the (analytic) natural
 projection. Let $\widetilde T$ denote the Nash bundle over
 $\widetilde Y$, i.e. the restriction over $\widetilde Y$ of the
 bundle $\bC^N \times U(d,N) \to \bC^N \times G(d,N)$, where $U(d,N)$
 is the tautological bundle over $G(d,N)$.  We consider a continuous,
 stratified vector field $\bv$ on a subset $V\subset Y$. The
 restriction of $\bv$ to $V$ has a well-defined canonical lifting
 $\tilde \bv$ to $\nu^{-1}(V)$ as a section of the Nash bundle
 $\widetilde T \to \widetilde Y$ (see e.g. \cite{BS}, Prop.9.1).

  We say that the stratified vector field $\bv$ on $Y$ is {\em radial-at-infinity} if it is defined on the restriction to $Y$ of the
  complement of a sufficiently large ball $B_M$ centered at the origin
  of $\bC^N$, and if $\bv$ is transversal to $S_R$, pointing outwards,
  for any $R>M$.  The ``sufficiently large'' radius $M$ is furnished
  by the following well-known result.
 
\begin{lemma}\label{l:link} 
There exists $M\in \bR$ such that, for any
$R\ge M$, the sphere $S_R$ centered at the origin of $\bC^N$ and of
radius $R$ is stratified transversal to $Y$, i.e. transversal to all
strata of the stratification $\cW$.  \qed
\end{lemma} This follows essentially from Milnor's finiteness result
 \cite[Cor. 2.8]{Mi} applied to the strata of $Y$. We shall call
$K_\ity(Y) := Y \cap S_R$ the {\em link at infinity} of $Y$.  By Lemma \ref{l:link} and standard isotopy arguments, $K_\ity(Y)$ does not depend on the radius
$R$, provided that $R>M$. Actually the link at infinity $K_\ity(Y)$
does not depend on the choice of the center of the sphere either,
since for two spheres $S', S''$ centered at two different points, the
links $Y \cap S'$ and $Y \cap S''$ are isotopic, provided the radii
are large enough.

An important consequence of Lemma
\ref{l:link} is that every vector field $\bw$ defined on  $Y \cap S_R$ without zeros, can be extended to the exterior of the ball $Y\cap
B_R$ without zeros, via the fibration provided by Lemma
\ref{l:link}. Therefore the obstruction to extend $\bw$ to the whole $Y$
subsists only inside the ball $Y\cap B_R$.
  This motivates the following definition, yet related to Dubson's definition
  of Euler obstruction \cite{Du1}:

\begin{definition} \label{d:euler} Let $\tilde \bv$ be the lifting to
a section of the Nash bundle $\widetilde T$ of a radial-at-infinity, stratified
vector field $\bv$ over $K_\ity(Y)= Y \cap S_R$. We call {\em global Euler
obstruction of} $Y$, and denote it by $\Eu(Y)$,
the obstruction for extending $\tilde \bv$ as a nowhere zero section
of $\widetilde T$ within $\nu^{-1}(Y\cap B_R)$.
\end{definition}
 
To be precise, the obstruction to extend $\tilde \bv$ as a nowhere
zero section of $\widetilde T$ within $\nu^{-1}(Y\cap B_R)$ is in fact
a relative cohomology class $o(\tilde \bv) \in H^{2d}(\nu^{-1}(Y \cap B_R),
\nu^{-1}(Y \cap S_R))$.  The Euler obstruction of $Y$  is the evaluation of ${ o} (\tilde \bv)$ on the fundamental
class of the pair $(\nu^{-1}(Y\cap B_R), \nu^{-1}(Y\cap S_R))$.  Thus $\Eu(Y)$ is an integer. By the preceding discussion, $\Eu(Y)$ does not depend on the
 radius of the sphere defining the link at infinity $K_\ity(Y)$.
 Since two radial vector fields are homotopic as stratified vector
 fields, it does not depend on the choice of $\bv$ either. Elementary
 obstruction theory tells us that $o(\tilde \bv)$ is also independent
 of the way we extend the section $\tilde \bv$ to $\nu^{-1}(Y\cap
 B_R)$, see \cite{St}. Moreover, as shown in \cite{MP} or \cite{Du1},
 this is also independent on the blow-up $\nu$ (one can work with any blow-up which extends the tangent bundle over $Y_\reg$). 
\subsection{Properties of $\Eu(Y)$}\label{ss:prop}
Let us recall that the local Euler obstruction is defined at
each point of $Y$, and it is constant on each stratum of a Whitney
stratification (see for instance \cite{Du1,BS,LT}). Thus, given a stratum
$W_i\subset Y$ of $\cW$, one denotes by $\Eu_Y(W_i)$ the local Euler obstruction
at some point of $W_i$. Our just defined global Euler obstruction has the following properties.
\begin{enumerate}
\item If $Y$ is non-singular, then $\Eu(Y)= \chi(Y)$, the Euler-Poincar\'e characteristic of $Y$. 
\item $\Eu(Y)= c_M^d(Y)$, the top degree Chern-Mather class of $Y$.
\item $ \Eu(Y) = \sum_{\cW_i\subset Y} \chi(\cW_i) \Eu_Y(\cW_i)$.
\end{enumerate}
The property (a) is clear from the definition. For property 
(b), we notice that our obstruction $o(\tilde \bv)$ in cohomology is
evaluated on the homology orientation class. This gives a class in $H_0(\tilde Y)$, so it is an integer. The
map $\nu_* : H_0(\tilde Y)\to H_0(Y)$ is obviously an
isomorphism and takes $\Eu(Y)$ into the top Chern-Mather
class $c_M^d(Y)$.\\ 
 To explain the equality (c), let us recall the relation to Dubson's definition of the Euler obstruction of some analytic variety relative to a non-characteristic open set. In the language of \cite{Du2}, our ball $B_R$
 is a non-characteristic open set with respect to the Whitney stratification,
 which just means that $S_R$ is transversal to strata. The property (c) tells that $\Eu(Y)$ is an Euler-Poincar\'e characteristic weighted by
 the constructible function $\Eu_Y$. 
  The equality follows immediately by Dubson's \cite[Theorem 1]{Du2} applied to our setting. In terms of constructible functions, the equality follows from a direct image argument and a complete proof in this spirit can be found in \cite[(5.65)]{Sch-book}. A different proof can be derived from \cite[Theorem 4.1]{BS}.  One may imagine still another variant of proof, based on extending vector fields, along the lines of the proof in \cite{BLS}, see our comments in \S \ref{s:remarks}.


\section{Affine Lefschetz pencils and Main Theorem}

 
 In the local setting, various authors
have proved ``Lefschetz type'' formulas for the local Euler obstruction, see for instance \cite{Du1,LT,BLS,Sch,BMPS}; we refer to the bibliography for background on this topic. 

We now introduce global affine Lefchetz pencils and we define the set
 of invariants that enter in our formula for $\Eu(Y)$.
  Let us assume that the coordinates of $\bC^N$ are fixed and recall that
 $H_\ity$ is the hyperplane at infinity $\bP^N \setminus \bC^N$.  An {\em
 affine pencil of hyperplanes} $H_t$, $t\in \bC$, is defined by a
 linear function $l: \bC^N \to \bC$, where $H_t := l^{-1}(t)$.  The intersection $A:= \overline{H_t} \cap H_\ity$ is the same for all $t$ and it is called {\em axis} of the
 pencil.
 We need to work with generic pencils and therefore recall below a well-known result on the existence of such pencils, see
e.g. \cite{La}, \cite{GM} or
\cite{Ti-ijm} for more general statements.  
\begin{lemma}\label{p:generic}
 There exists a Zariski open dense
subset $\Omega$ of linear functions on $\bC^N$ such that, for any
$l\in \Omega$: 
\begin{enumerate} 
\item the axis $A$ of the pencil
defined by $l$ is transversal to all the strata of $\cW$ which are
contained in the hyperplane at infinity $H_\ity$. 
\item there exists
a finite set $B\subset \bC$ such that, for all $t\in \bC \setminus B$,
the hyperplane $H_t$ cuts transversally all the strata of $Y$. For
$t\in B$, the hyperplane $H_t$ is transversal to the strata of $Y$ at
all points except finitely many, which are stratified Morse
singularities of the function $l$.\qed \end{enumerate}
\end{lemma}
\begin{definition}\label{d:pencil} The pencil of hyperplanes defined
by $l\in \Omega$ will be called {\em a Lefschetz pencil with respect to $Y$}.
\end{definition}
Let us also recall
the definition of complex stratified Morse singularities, used in Lemma \ref{p:generic}, since this is important in the proof of our main result.


\begin{definition}\label{d:morse} (Lazzeri '73, Benedetti '77, Pignoni
'79, Goresky-MacPherson '83 \cite[p.52]{GM}.) Let $\cW$ be a local Whitney stratification of a germ $(X,x_0)\subset (\bC^N,x_0)$ of a complex 
analytic space. Let $f: (X,
x_0) \to \bC$ be a holomorphic function germ and let $F: (\bC^N,
x_0) \to \bC$ denote some extension of it. We say that $f$ is a {\em general function} at $x_0$ if $\d F_{x_0}$ does not
vanish on any limit of tangent spaces to $W_i$, $\forall i\not= 0$,
and to $W_0\setminus \{ x_0\}$, where $W_0$ denotes the stratum to which $x_0$ belongs.  One says that $f: (X, x_0) \to \bC$
is a {\em stratified Morse function germ} if: $\dim W_0 \ge 1$, $f$ is
general with respect to the strata $W_i$, $i\not= 0$ and the
restriction $f_{|_ {W_0}}: W_0 \to \bC$ has a Morse point at $x_0$.
\end{definition}


\subsection{The global polar invariants}\label{ss:polar}

It appears that the polar invariants play a key role in studying the topology of spaces, sometimes in connection to the Lefschetz slicing method. Local polar multiplicities enter in the 
local Euler obstruction formula of L\^e-Teissier \cite[(5.1.2)]{LT}.
Polar classes (cf \cite{LT, Pi}) and Chern-Mather classes determine each other via a Todd type formula proved by Ragni Piene \cite{Pi}.
 
 We now define the global invariants which enter in our formula.
 The first number $\alpha_Y^{(1)}$ is by definition the number of Morse points on the regular part $Y_\reg$ of a Lefschetz pencil on $Y$.
 Next step, consider a general hyperplane $H_t\cap Y$ of the Lefschetz pencil and take a new Lefschetz pencil of $H_t\cap Y$: we get the second number 
$\alpha_Y^{(2)} := \#$ Morse points of the second Lefschetz pencil on $H_t\cap Y_\reg$.
This continues by induction and we get a sequence of non-negative integers:
\[ \alpha_Y^{(1)}, \alpha_Y^{(2)},\cdots , \alpha_Y^{(d)},\]
to which we attach the last one $\alpha_Y^{(d+1)}:=$ the number of points
of the intersection of $Y_\reg$ with a generic codimension $d$ plane in $\bC^N$. This is in fact just the degree of $Y$.

All these numbers are well-defined invariants of $Y$, by the connectivity of the Zariski open sets of generic slices and of pencils which we use.
In fact, these invariants can be interpreted as global polar multiplicities,
similar to the local ones used by L\^e-Teissier \cite{LT}, see also Piene \cite{Pi}.
 
 Global polar invariants have been also used by the second named author in order to characterize a certain equisingularity at infinity of families of affine hypersurfaces \cite{Ti-equi}. 
 
 Our result can be stated as follows:
\begin{theorem}\label{t:main}
If $Y\subset \bC^N$ is an algebraic variety  of pure dimension $d$, then
its global Euler obstruction is $\Eu(Y) = \sum_{i=1}^{d+1} (-1)^{d-i+1} \alpha_Y^{(i)}$.

\end{theorem} 
 
 Let us stress that only Morse points on the regular part (or the regular part of repeated slices) come into the description of the Euler obstruction $\Eu(Y)$, in other words the Morse stratified singularities on 
 other strata are ignored. Consequently, the stratification $\cW$ does not appear in
 the statement at all.
 Furthermore, the formula resembles to an Euler characteristic formula,
 by attaching cells via Lefschetz slicing method. In general, it is not the Euler characteristic of $Y$ or of $Y_\reg$. Nevertheless, in case $Y$ is  non-singular, i.e. if $Y = Y_\reg$, we clearly have $\Eu(Y) = \chi(Y)$.


\section{Proof of the Main Theorem}\label{proof}
 
 Let us fix a Lefschetz pencil with respect to $Y$, of axis $A$, given by a linear
 function $l\in\Omega$ (Lemma
 \ref{p:generic}, Definition \ref{d:pencil}). Its restriction to $Y$ is a stratified submersion
 away from a finite set of points $\Sigma = \{ y_1, \ldots , y_k\}$.
 At each such point the function germ $l_{|Y} : (Y, y_i) \to \bC$ has
 a stratified Morse singularity. Let us then take a large enough disk
 $D\subset \bC$ such that the finite set $B$ of critical values of $l_{|Y}$
 is contained in $D$.  


 We fix a generic slice of $H_t := l^{-1}(t)$, i.e. $t\in D\setminus B$.
 Notice that the intersections of $H_t$ with the strata of $\cW$ constitutes a Whitney stratification of $Y\cap H_t$, to which we refer in the following. We apply Lemma \ref{l:link} to $Y$ and then to $Y\cap H_t$. This yields 
 a positive real $M$ such that the sphere $S_R$ is stratified transversal to $Y$ and to $Y\cap H_t$, for all $R\ge M$. We then fix an $R\ge M$.
 
  We remark that the critical set $\Sigma$ is contained
 in $Y\cap B_R \cap l^{-1}(D)$.  We shall construct a special continuous
 stratified vector field $\bv$ on $Y\cap B_R \cap l^{-1}(D)$ which is
 radial on $Y\cap S_R \cap l^{-1}(D)$ and points outwards the tube
 $l^{-1}(\partial D)$. We notice that the procedure used
 in the local case in \cite{BLS,BMPS} applies to our situation,
 when replacing a small Milnor ball by our big ball $B_R$. Our function
 will be the global function $l$, instead of a germ, and this
 function has several stratified critical points.

\begin{figure}[hbtp]\label{f:vectorfield}
\begin{center}
\epsfxsize=5cm
\leavevmode
\epsffile{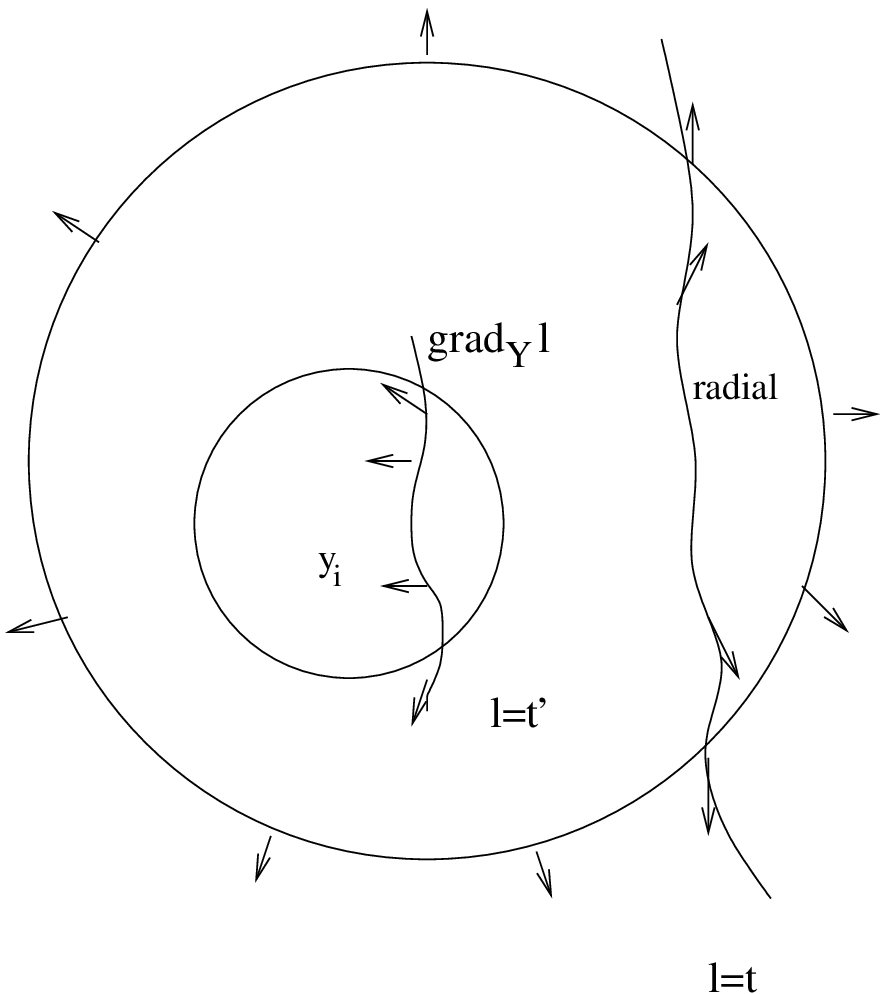}
\end{center}
\caption{} 
\label{f:1}   
\end{figure}

 As a matter of fact, we start with a radial-at-infinity vector field $\bv$ on $Y\cap S_R \cap l^{-1}(D)$ which in addition
 is radial-at-infinity with respect to the strata of the
 fixed slice $Y\cap H_t$. Next we extend this vector field 
 over $Y\cap H_t\cap B_R$ (one can always do this such that
 the extension has isolated zeros).  Then, as shown in
 {\em loc.cit.}, one may extend this to a continuous stratified vector field without zeroes, outside $Y\cap H_t\cap B_R$ and within a tubular neighbourhood of $Y\cap H_t\cap B_R$, by plugging in the lifting by $l$ of a vector
 field on the disk $D$ which is radial from the point $t$. We may
 further extend and deform this without zeros outside the tubular neighbourhood such that it becomes the
 gradient vector field $\grad_Y l$ (defined below) in the neighbourhood of any point $y_i \in \Sigma$.   See Figure \ref{f:1}.

Let us first give the precise definition of $\grad_Y l$. 
Take the {\em complex conjugate} of the gradient of $l$ and project
 it to the tangent spaces of the strata of $Y$ into a stratified vector field on $Y$. This 
  may be not continuous, but one can make it continuous by ``tempering'' it in the neighbourhood of ``smaller'' strata. One then gets a 
continuous stratified vector field, well-defined up to stratified homotopy, which
is by definition $\grad_Y l$ and which we have called above {\it the gradient
vector field}. 

Note that the only zeros of $\grad_Y l$ occur precisely at the points $y_i\in \Sigma$.
   If $\nu: \widetilde Y \to Y$ is the Nash
 blow-up of $Y$ and $\widetilde T$ is the Nash bundle over $\widetilde
 Y$, then $\grad_Y l$ lifts canonically to a never zero section
 $\widetilde{\grad_Y l}$ of $\widetilde T$ restricted to $\nu^{-1}(Y \cap S_{\varepsilon})$, where $S_{\varepsilon}$ is a
small Milnor sphere around some fixed $y_i$. In our problem, we are interested in the obstruction to extend
$\widetilde{\grad_Y l}$ without zeros throughout $\nu^{-1}(Y \cap B_{\varepsilon})$. One considers a more general problem in \cite{BMPS}:  given a holomorphic function germ $f: (X, x_0) \to \bC$, compute the obstruction of extending $\widetilde{\grad_X f}$ without zeros from $\nu^{-1}(X \cap S_{\varepsilon})$ throughout $\nu^{-1}(X \cap B_{\varepsilon})$. This is called {\em local Euler obstruction of $f$}, or ``defect of $f$'', and is denoted by $\Eu_f(X,x_0)$. In case of Morse singular points, one may directly compute this obstruction as follows. 
\begin{lemma}\label{l:morse} 
Let $l$ be a holomorphic function germ
on $(Y, y)$ with a stratified Morse singularity at $y \in W_0$.
Then the local Euler obstruction of $l$ is $0$ if $\dim W_0 < \dim Y$, and
is $(-1)^{\dim_\bC Y}$ if $W_0 = Y_\reg$.
\end{lemma}

\begin{proof} Take a small enough ball $B_\e$ in $Y\subset \bC^N$,
centered at $y$ and of radius $\e >0$. Let $v$ be the gradient
vector field $\grad_Y l$ restricted to the sphere $Y \cap \partial B_\e$
and consider the tautological lift $\tilde v$ to the Nash blow-up
$\widetilde Y$. From the definition of a Morse function we have that
the form $\d l$ does not vanish on any limit of tangent spaces at
points in the regular stratum $Y_\reg$.  Since $\widetilde Y$ is
obtained by attaching to $TY_\reg$ all limits of tangent spaces of
points in $Y_\reg$, one has that if $y \notin Y_\reg$, then, by the Definition \ref{d:morse} of stratified Morse points,  the
section $\tilde v$ of $\widetilde T$ can be extended over $\nu^{-1}(Y
\cap B_\e)$ without zeros. This proves the first part of our statement, since the extension can be done as follows. Let us think of ${\widetilde T}$ as being a subset
of $({\bC^N} \times G(d,N)) \times
\bC^N$. At each point $(x, H)\in \widetilde Y$ one adds to $\tilde
v$ the tautological lift of the vector 
\[ \gamma(\| x\|) \cdot \proj_H (\overline{\grad l}),\]
 where $\proj_H$ denotes the projection to $H$
and $\gamma$ denotes a continuous non-negative real function defined
on $[0, \e]$ with values in $[0,1]$, such that $\gamma(\e)=0$ and
$\gamma(0)=1$ (for instance $\gamma$ can be taken linear.)

The proof of the second part of our statement, for $y \in Y_\reg$, goes as follows. Locally, the Nash bundle is the usual
tangent bundle of $Y_\reg$ and $\Eu_l(Y,y)$ is by definition the
Poincar\'e-Hopf index of $grad_Y l$ at $y$. One deduces from \cite[Th.7.2]{Mi} that $\Eu_l(Y,y) =(-1)^{\dim_\bC Y}$, since $y$ is just a complex Morse point in the classical sense. 
\end{proof}
The above lemma is natural since the
Euler obstruction is defined via the Nash blow-up and the latter only
takes into account the closure of the tangent bundle over the regular
part $Y_\reg$.  

We pursue the proof of Theorem \ref{t:main}.
 We have shown how one may extend the radial-at-infinity vector field $\bv$ to a vector field over $Y\cap B_R \cap l^{-1}(D)$. In the same manner as in \cite{BLS,BMPS} for the local case, one may extend this further without zeros on
$(Y\cap B_R) \setminus (Y\cap B_R \cap l^{-1}(D))$. 

 The construction we have done shows, firstly, that the global Euler obstruction
$\Eu(Y)$ is exactly the obstruction
to extend $\tilde \bv$ inside the bounded tube $\nu^{-1}(Y\cap B_R
\cap l^{-1}(D))$ and secondly, that this obstruction is
 the sum of two terms:
 
\begin{enumerate}
\item the obstruction to extend $\tilde \bv$ within the slice $\nu^{-1}(Y\cap H_t \cap B_R)$, as a lift of a stratified vector field with respect to the strata of $Y\cap H_t$.

\item the obstructions due to the isolated zeros of the gradient vector field
 $\grad_Y l$. 
 
\end{enumerate}
The obstruction (a) is $\Eu(Y\cap H_t)$ essentially by definition, since
the Whitney strata of $Y\cap H_t$ are precisely the intersections of the
strata of $Y$ with $H_t$, by the assumed transversality of $H_t$ to the
Whitney stratification $\cW$ of $Y$.
 For (b), by Lemma \ref{l:morse}, the local obstruction at some
 stratified Morse point $y_i$ is zero if $y_i$ is not on $Y_\reg$ and it is
 $(-1)^d$ for $y_i\in Y_\reg$. So the sum of all such local obstructions
 is equal to $(-1)^{d}$ times the number of Morse points.
 This number is just the global invariant $\alpha_Y^{(1)}$ defined at \ref{ss:polar}.  We therefore get:
  \begin{equation} \label{eq:lef}
  \Eu(Y) = \Eu(Y\cap H_t) + (-1)^{d}\alpha_Y^{(1)}
\end{equation}

  We may now apply the same procedure to $Y\cap H_t$ instead of $Y$, get a formula similar to (\ref{eq:lef}) for $Y\cap H_t$, i.e. $\Eu(Y\cap H_t) = \Eu(Y\cap H_t \cap H_{t'}) + (-1)^{d-1}\alpha_Y^{(2)}$,
  and so on by induction. Altogether this leads to formula (\ref{eq:main}). The proof of Theorem \ref{t:main} is now complete.
 

\section{Further remarks}\label{s:remarks}


We have seen that our formula (\ref{eq:main}), proved via (\ref{eq:lef}), belongs to the vein of Lefschetz type formulas and in the same time it is a L\^e-Teissier type one. Let us see how the term $\Eu(Y\cap H_t)$ in formula (\ref{eq:lef}) can be refined in another direction. In the local setting, the obstruction to extend $\tilde \bv$ within the slice $\nu^{-1}(X\cap H_t \cap B_\e)$, as a lift of a stratified vector field with respect to the strata of $X\cap H_t$, has been looked up in \cite{BLS,BMPS}. The proof of the main formula in \cite{BLS} can be rephrased as follows: one first shows that actually $\Eu_X(x_0) = \Eu(X\cap B_\e \cap H_t)$ and next that $\Eu(X\cap B_\e \cap H_t)$ can be expressed as a weighted Euler characteristic. In the global setting, the latter fact has the following interpretation:
\begin{equation}\label{eq:slice}
\Eu(Y\cap H_t) = \sum_{W_i\subset Y}
\chi(Y\cap B_R\cap H_t\cap W_i) \cdot \Eu_Y(W_i),
\end{equation} 
By the choice of the radius $R$ in the beginning of the proof of our Theorem \ref{t:main},  we have that $Y\cap B_R\cap H_t\cap W_i$ 
is diffeomorphic to $Y\cap H_t\cap W_i$, hence the sum (\ref{eq:slice}) is in turn equal to
$\sum_{W_i\subset Y} \chi(Y\cap H_t\cap W_i) \cdot \Eu_Y(W_i)$.
 Formula (\ref{eq:lef}) becomes therefore:
\begin{equation}\label{eq:eu} 
\Eu(Y) = \sum_{W_i\subset Y}
[ \chi(H_t\cap W_i) \cdot \Eu_Y(W_i) ] + (-1)^{d} \alpha_Y^{(1)}. 
\end{equation} 
 One may alternatively prove formula (\ref{eq:slice}) by using property (c) in \S \ref{ss:prop}.
 
 Actually, Dubson's definition of Euler obstruction relative to a non-characteristic real analytic open set \cite{Du1, Du2} together with 
 the Lefschetz slicing method allows one to understand in a unitary way the local or global formulas for the Euler obstruction.
For instance, the formula proved in \cite{BLS} could be rephrased saying that the local Euler obstruction $\Eu_X(x_0)$ equals the  
Euler obstruction of the complex link of $X$ at $x_0$ to which one applies  property (c) in \S \ref{ss:prop}.
 
 Alternatively, instead of applying property (c), one may use local Lefschetz slicing in the similar manner as in the global setting (proof of Theorem \ref{t:main}). By applying formula (\ref{eq:lef}) repeatedly to the  
Euler obstruction of the complex link of $X$ at $x_0$, we get exactly L\^e-Teissier's formula for the local Euler obstruction \cite{LT}.


\begin{thebibliography}{MMM}



\bibitem[BS]{BS}
J.P. Brasselet and M.H. Schwartz {\em Sur les classes de Chern d'un
ensemble analytique complexe} Asterisque {\bf 82-83} (1981), 93-147.

\bibitem[BLS]{BLS}
J.-P. Brasselet, L\^e D.T., J. Seade, {\em Euler obstruction and indices
of vector fields}, Topology 39, no. 6 (2000), 1193--1208.


\bibitem[BMPS]{BMPS}
J.-P. Brasselet, D.B. Massey, A.J. Parameswaran, J. Seade, {\em Euler
obstruction and defects of functions on singular varieties},
math.AG/0302238, to appear in J. London Math. Soc.


\bibitem[Du1]{Du1} 
A. Dubson, {\em Classes caract\'eristiques des vari\'et\'es
singuli\'eres}, C.R. Acad.Sc. Paris {\bf 287} (1978), no. 4,  237--240.

\bibitem[Du2]{Du2}
 A. Dubson, {\em Formule pour l'indice des complexes constructibles et
 des D-modules holonomes}, C. R. Acad. Sci. Paris S\'er. I Math. 298
 (1984), no. 6, 113--116.

\bibitem[GM]{GM}
M. Goresky, R. MacPherson, {\em Stratified Morse theory}, Ergebnisse
der Mathematik und ihrer Grenzgebiete. 3. Folge, Bd. 14.  Berlin
Springer-Verlag 1988.


\bibitem[La]{La}
 K. Lamotke, {\em The topology of complex projective varieties after
 S. Lefschetz}, Topology 20 (1981), 15-51.


\bibitem[LT]{LT} L\^ e D.T., B. Teissier {\em Vari\'et\'es polaires
locales et classes de Chern des vari\'et\'es singulieres}, Ann. of
Math {\bf 114} (1981), 457-491.


\bibitem[MP]{MP} R. MacPherson, {\em  Chern classes for singular varieties},
Ann. of Math {\bf 100} (1974), 423-432.


\bibitem[Mi]{Mi}
J. Milnor, {\em Singular points of complex hypersurfaces}, Ann. of
Math. Studies~61, Princeton 1968.

\bibitem[Pi]{Pi}
 R. Piene, {\em Polar classes of singular varieties}, Ann. Sci. Ecole Norm. Sup. (4) 11 (1978), no. 2, 247--276.


\bibitem[Sch1]{Sch}
J. Sch\"urmann, {\em A short proof of a formula of Brasselet, L\^e and
Seade}, math. AG/ 0201316.

\bibitem[Sch2]{Sch-book}
J. Sch\"urmann, {\em Topology of singular spaces and constructible sheaves}, Mathematical Monographs-new series, Birkh\" auser 2003. 

\bibitem[St]{St}
N. Steenrod, {\em The Topology of Fiber Bundles}, Princeton
Univ. Press, 1951.

 
\bibitem[Ti1]{Ti-equi}
M. Tib\u ar, {\em Asymptotic equisingularity and topology of complex hypersurfaces}, Internat. Math. Res. Notices 1998, no. 18,
  979--990. 
  
 \bibitem[Ti2]{Ti-ijm}
M. Tib\u{a}r, {\em Connectivity via nongeneric pencils},
Internat. J. Math. 13 (2002), no. 2, 111--123.

 
\end{thebibliography}
\end{document}